\documentclass[12pt]{amsart}

\usepackage{latexsym}
\usepackage{amssymb}
\usepackage{amsmath}
\usepackage{amscd}
\usepackage{amsthm}
\usepackage[all]{xy}  

\numberwithin{equation}{section}

\newtheoremstyle{fancy1}{10pt}{10pt}{\itshape}{12pt}{\textsc\bgroup}{.\egroup}{8pt}{
}
\newtheoremstyle{fancy2}{10pt}{10pt}{}{12pt}{\itshape}{.}{8pt}{ }

\theoremstyle{fancy1}

\newtheorem*{thethm*}{Theorem A}

\newtheorem{main}{Theorem}
\newtheorem*{main*}{Result}

\newtheorem*{cor*}{Corollary}

\theoremstyle{fancy2}

\newtheorem*{def*}{Definition}

\newtheorem*{rem*}{Remark}
\newtheorem*{rems*}{Remarks}

\newtheorem*{example*}{Example}

\newtheorem*{examples*}{Examples}

\theoremstyle{remark}

\newcommand{\cref}[1]{Corollary~\ref{#1}}

\newcommand{\CP}{\mathbb{C\mkern1mu P}}

\newcommand{\C}{{\mathbb{C}}}
\newcommand{\R}{{\mathbb{R}}}
\newcommand{\Z}{{\mathbb{Z}}}
\renewcommand{\H}{{\mathbb{H}}}

\newcommand{\SO}{\ensuremath{\operatorname{\mathsf{SO}}}}

\newcommand{\U}{\ensuremath{\operatorname{\mathsf{U}}}}
\newcommand{\SU}{\ensuremath{\operatorname{\mathsf{SU}}}}

\renewcommand{\S}{\ensuremath{\operatorname{\mathsf{S}}}}

\def\con#1=#2(#3){#1 \equiv #2 \bmod{#3}}

\newcommand{\set}[1]{\left\{#1\right\}}
\newcommand{\grp}[1]{\langle #1 \rangle}
\newcommand{\ab}{\text{ab}}  

\newcommand{\diag}{\ensuremath{\operatorname{diag}}}

\newcommand{\del}{\partial}

 \DeclareMathOperator{\pr}{pr}

\DeclareMathOperator{\Id}{Id}

\newcommand{\bd}{\overline{d}}

\newcommand{\bwp}{\overline{\wp}}

\newcommand{\bE}{\overline{E}}

\begin{document}


\title[6-Dimensional Cohomogeneity One Manifolds]{Diffeomorphism type of 6-dimensional\\cohomogeneity one manifolds}

\author{Corey A. Hoelscher}
\address{Rutgers University\\
New Brunswick, NJ}

\begin{abstract}
In this paper we give a full diffeomorphism characterization of compact simply connected cohomogeneity one manifolds in dimension six.
\end{abstract}

\maketitle

%
%
%
%
%
%

\renewcommand{\thetable}{\Roman{table}}

\section*{Introduction}

A cohomogeneity one manifold is a smooth manifold with a large symmetry group, in the sense that there is a compact Lie group that acts on the manifold with a one dimensional orbit space, or equivalently with at least one orbit of codimension one. Since the orbits of any smooth compact Lie group action are embedded submanifolds of some codimension, cohomogeneity one manifolds are the manifolds with the next ``largest'' symmetry group, after homogeneous spaces, in this sense.

Cohomogeneity one manifolds have a natural importance in the theory of group actions on manifolds, and they are also important in many other areas of geometry. For example in Riemannian geometry, it was shown in \cite{GZ2000} that many cohomogeneity one manifolds admit metrics of non-negative sectional curvature, including every principal $S^3$ bundle over $S^4$. Cohomogeneity one manifolds also give examples of positive and non-negative Ricci curvature \cite{GZ2002} and recently a metric of positive sectional curvature was constructed on a manifold homeomorphic to the unit tangent bundle of $S^4$, using a cohomogeneity one structure. Cohomogeneity one manifolds are also important in other areas of geometry and mathematical physics. For example, they give new examples of Einstein and Einstein-Sasaki manifolds \cite{Conti,GHY} and examples of manifolds with exceptional holonomy groups \cite{CS,CGLP}.

Cohomogeneity one manifolds were classified in dimensions 4 and lower in \cite{Parker} and \cite{Ne}, and later in \cite{thesis} the compact simply connected ones were classified in dimensions 5, 6 and 7. These classifications describe the possible manifolds by expressing them as the union of two disk bundles glued along their common boundary. However it is not always clear which manifolds result from this construction. In dimensions 4 and lower, it is known which compact simply connected manifolds result \cite{Parker,Ne} and in dimension 5 a diffeomorphism characterization is given in \cite{thesis}. The main result of this paper settles this question in dimension 6.


\begin{main}\label{6dimthrm}
A compact simply connected 6-dimensional cohomogeneity one manifold is diffeomorphic to one of the following: 
\begin{enumerate}
	\item[i.] a compact symmetric space;
	\item[ii.] the nontrivial $S^4$ sphere bundle over $S^2$;
	\item[iii.] a nontrivial $\CP^2$ bundle over $S^2$ with structure group $\SU(3)/\Z_3$; or
	\item[iv.] an $S^2$ sphere bundle over $\CP^2$ or $S^2\times S^2$ whose structure group reduces to $\SO(2)$.
\end{enumerate}
Conversely all such manifolds admit cohomogeneity one actions.
\end{main}

An interesting consequence of the proof of Theorem \ref{6dimthrm} and of \cite{thesis} is that among the 6-dimensional compact simply connected symmetric spaces, $S^6$, $\CP^3$, $\SO(5)/\SO(2)\SO(3)$ and $S^3\times S^3$ admit only isometric cohomogeneity one actions whereas $S^4\times S^2$, $\CP^2\times S^2$ and $S^2\times S^2\times S^2$ admit both isometric and non-isometric cohomogeneity one actions.

A characterization similar to Theorem \ref{6dimthrm} in dimension 7 would be much more difficult. The first step in this direction is taken in \cite{cohom1_top} where the homology groups of these manifolds are computed. 

The paper is organized as follows. In the first section, we review the general structure of cohomogeneity one manifolds and we recall a certain fiber bundle structure associated to non-primitive actions. In Section 2, we consider each remaining family of manifolds from the classification and use this non-primitivity fiber bundle structure to determine their diffeomorphism types. A summary of this information is given in Table \ref{t:summary} for the convenience of the reader.


The author is happy to express thanks to W.~Ziller and S.~Ferry for several helpful discussions and to the former for suggesting improvements to this manuscript.

%
%
%
%
%
%

\renewcommand{\thetable}{\theequation}

\section{Cohomogeneity One Manifolds}\label{prelims}

In this section we recall the basic structure of general cohomogeneity one manifolds, and the bundle structure of the non-primitive ones. This bundle structure is the key to proving Theorem \ref{6dimthrm}.

\subsection{Basic Structure}
Suppose $G$ is a compact connected Lie group which acts on the compact connected manifold $M$ by cohomogeneity one. $M$ is called an \emph{interval cohomogeneity one manifold} if the orbit space $M/G$ is a closed interval. This happens whenever $M$ has finite fundamental group, for example. Such a manifold $M$ is determined by its group diagram, $G\supset K^-,K^+\supset H$, where $H$ is a principal isotropy subgroup and $K^\pm$ are certain corresponding non-principal isotropy subgroups with the property $K^\pm/H\simeq S^{\ell_\pm}$. (See \cite{thesis}, \cite{GZ2000} or \cite{GWZ} for more details.) More precisely $M$ can be decomposed $G$-equivariantly as
\begin{equation}\label{decomp}
  M\, = \, G\times_{K^-}D_- \,\, \cup \,\, G\times_{K^+}D_+ \text{ \hspace{1em} where \hspace{1em} } S^{\ell_\pm}=\del D_\pm \simeq K^\pm/H
\end{equation}
and where the two halves $G\times_{K^\pm}D_\pm$ are glued along their common boundary $\del(G\times_{K^\pm}D_\pm)=G\times_{K^\pm} K^\pm/H\simeq G/H$. Conversely any collection of compact groups $G\supset K^-,K^+\supset H$ with $K^\pm/H\simeq S^{\ell_\pm}$ determines an interval cohomogeneity one manifold given by \eqref{decomp}.

\subsection{Non-primitive actions}\label{sec:nonprim}


Recall that a cohomogeneity one manifold $M$ is called non-primitive if it has some group diagram $G\supset K^-,K^+\supset H$ such that there is a closed connected proper subgroup $L\subset G$ which contains $K^-$, $K^+$ and $H$. Then $L\supset K^-,K^+\supset H$ is a valid group diagram corresponding to some cohomogeneity one manifold $M_L$. Further $M$ is $G$-equivariantly diffeomorphic to $G\times_L M_L$ and in particular we have the \emph{non-primitivity fiber bundle}
\begin{equation}\label{nonprim}
M_L \to M \to G/L.
\end{equation}
The structure group for this bundle is $L$, though not necessarily effectively. See \cite{thesis} for details.

We can potentially enlarge the effective structure group as follows. Suppose there is a compact connected Lie group $J$ which acts effectively on $M_L$, and a Lie homomorphism $\phi: L\to J$, such that the $L$ and the $J$ actions on $M_L$ are compatible, i.e.~for $\ell\in L, x\in M_L$ we have $\ell\star x= \phi(\ell)\star x$. Define $P=G\times_L J=G\times J/\sim$ where $(g,j)\sim (g\ell^{-1},\phi(\ell)j)$ for $\ell\in L$ and note that $G$ acts on the left and $J$ on the right of $P$. We then see that $M$ is $G$-equivariantly diffeomorphic to $P\times_J M_L$, as in \cite[Prop. 1.5]{thesis}. Furthermore $M_L \to M\to G/L$ is the bundle associated to the principal $J$ bundle $P$ and the action of $J$ on $M_L$.

%
%
%
%
%
%

\section{Proof of Theorem \ref{6dimthrm}}\label{body}

In this section we will prove the main theorem from the introduction. Suppose that $M$ is a compact simply connected cohomogeneity one manifold of dimension 6. We will show that $M$ is diffeomorphic to one of the manifolds listed in Theorem \ref{6dimthrm}.

The classification of these manifolds \cite[Thm. A]{thesis} says that $M$ is diffeomorphic to one of the following: a symmetric space; the product of a homogeneous space and a cohomogeneity one manifold; or a manifold given by one of the group diagrams listed in Table II of \cite{thesis}. 

The first case is clear so suppose that $M$ is the product of a homogeneous space and a cohomogeneity one manifold. Since $M$ is simply connected and 6-dimensional, each factor in the product must have dimension 2, 3, or 4. The only compact simply connected homogeneous spaces in these dimensions are $S^2$, $S^3$, $S^4$, $S^2\times S^2$ and $\CP^2$. 
Furthermore, the only additional compact simply connected cohomogeneity one manifold in these dimensions is $\CP^2\# -\CP^2$ (\cite{Parker} and \cite{Ne}). Therefore the only possibility for $M$ that is not already a symmetric space is $S^2\times (\CP^2\# -\CP^2)$. Because $\CP^2\# -\CP^2$ is an $S^2$ bundle over $S^2$ with structure group $\SO(2)$ we see that $S^2\times (\CP^2\# -\CP^2)$ is an $S^2$ bundle over $S^2\times S^2$ with structure group $\SO(2)$. In particular this manifold already appears in Theorem \ref{6dimthrm}. 

In the rest of this section we consider the remaining cases: the manifolds given by the group diagrams from Table II of \cite{thesis}. For each of these families of manifolds we begin by recalling the group diagrams which describe the manifolds and the conditions on the group diagrams which ensure that the resulting manifolds are simply connected. All of these manifolds are non-primitive and in most cases we take advantage of the non-primitivity fiber bundle \eqref{nonprim} to determine their diffeomorphism type.

\subsection{Actions of type $N^6_A$:}\label{N^6_A}

  $$S^3\times S^1\times S^1 \,\, \supset \,\, \set{(e^{i a_-\theta},e^{ib_-\theta},e^{ic_-\theta})}\cdot H_+, \,\, \set{(e^{i a_+\theta},e^{ib_+\theta},e^{ic_+\theta})}\cdot H_- \,\, \supset \,\, H_-\cdot H_+$$
where $K^-\ne K^+$, $H_\pm=H\cap K^\pm_0$ is finite cyclic, $\gcd(b_\pm,c_\pm)=1$, $a_\pm=rb_\pm+sc_\pm$ for some $r,s\in \Z$, and $K^-_0\cap K^+_0\subset H$.

\vspace{.5em}

We show here that these are all isometric actions on $S^3\times S^3$. To see this, we first claim that $H=K^-\cap K^+$. Obviously $H\subset K^-\cap K^+$. Conversely, it is not hard to see that $K^-\cap K^+_0\subset H$ since $K^-=K^-_0\cdot(K^+_0\cap H)$ and $K^-_0\cap K^+_0\subset H$ by assumption. Similarly $K^-_0\cap K^+\subset H$. Now notice $L=K^-_0\cdot K^+_0=K^-\cdot K^+\simeq T^2$ since $K^\pm=K^\pm_0\cdot H_\mp$. Choose $l\in L$ and write $l=k_-k_+$ for $k_\pm\in K^\pm_0$. If $l\in K^-\cap K^+$ it is clear $k_-\in K^-_0\cap K^+\subset H$ and similarly $k_+\in H$, so $l\in H$. This proves $H=K^-\cap K^+$.

Now, \cite{thesis} describes the family of actions of $S^3\times S^1\times S^1$ on $S^3\times S^3$ given by $(g, z,w)\star (x,y) = (gx\bar z^r\bar w^s, \big((z^{c_-}\bar w^{b_-})^{n_-},(z^{c_+}\bar w^{b_+})^{n_+}\big)\star_1 y)$ where $\star_1$ is the usual torus action on $S^3$, and where we consider $S^1\subset S^3\subset \H$. The group diagram for this action is given by
$$S^3\times T^2  \supset  \set{(z^rw^s,z,w)|(z^{c_-}\bar w^{b_-})^{n_-}=1},  \set{(z^rw^s,z,w)|(z^{c_+}\bar w^{b_+})^{n_+}=1} \supset$$ $$\nonumber  \set{(z^rw^s,z,w)|(z^{c_-}\bar w^{b_-})^{n_-}\!=\!1\!=\!(z^{c_+}\bar w^{b_+})^{n_+}}.$$

It is clear that $\set{(z^rw^s,z,w)|(z^{c}\bar w^{b})^{n}=1}=\set{(z^a,z^b,z^c)|z\in S^1}\cdot\Z_n$ where $a=rb+sc$. Hence any $K^\pm$ from an action of type $N^6_A$ can be obtained in this way. Since $H=K^-\cap K^+$ it follows that any diagram of type $N^6_A$ is of this form. In particular, the resulting manifold is diffeomorphic to $S^3\times S^3$.

\subsection{Actions of type $N^6_B$:}

$$S^3\times S^3 \,\, \supset \,\, \set{(e^{i\theta},e^{i\phi})}, \,\, \set{(e^{i\theta},e^{i\phi})} \,\, \supset \,\, \set{(e^{ip\theta},e^{iq\theta})} \cdot \Z_n$$
where $\gcd(p,q)=1$ and $\Z_n\cap\set{(e^{ip\theta},e^{iq\theta})}=1$.

\vspace{.5em}

Let $L=K^\pm=T^2$. Then $G/L\simeq S^2\times S^2$ and $M_L\simeq S^2$ where $L$ acts on $M_L$ via $\phi:L\to \SO(3):(e^{i\alpha},e^{i\beta})\mapsto R\big(n(-q\alpha+p\beta)\big)$. Here $R(\theta)\in \SO(3)$ denotes rotation by $\theta$ about some fixed axis. Therefore the non-primitivity bundle for $M$ is $S^2\to M\to S^2\times S^2$ and the effective structure group can be taken to be $\SO(2)=\set{R(\theta)}\subset \SO(3)$.

Consider the principal $\SO(2)$ bundle $\wp:P\to S^2\times S^2$ associated to the non-primitivity sphere bundle $S^2\to M\to S^2\times S^2$. We know $P=G\times_L\SO(2)$ from Section \ref{sec:nonprim}. It is well known that such principal $\SO(2)$ bundles are characterized by their Euler class $e$ which we will now explicitly compute up to sign.

Notice that $P= G\times \SO(2)/F$ where $F=\set{(\ell^{-1},\phi(\ell))|\ell\in L}$ since $\SO(2)$ is abelian. Furthermore, the map $\wp:P\to G/L$ can be expressed as the projection $\wp:G\times \SO(2)/F \to G\times \SO(2)/T$, where $T=L\times \SO(2)$, since $G\times \SO(2)/T\simeq G/L$ naturally. Now consider the following diagram of bundles
\begin{equation}\label{eulerdiagram}
\displaystyle{\begin{array}{ccccc}
      T & \to & G\times \SO(2) & \to & G\times \SO(2)/T\\
      \downarrow & & \downarrow & & \downarrow \\
      T/F & \to & G\times \SO(2)/F & \to & G\times \SO(2)/T
 \end{array}  }
\end{equation}
where the bottom bundle is $\wp$. Denote the top bundle by $\bwp$ and denote the vertical maps by $\pi$. It is clear from the Gysin sequence for $\wp$ that $\ker\big(\wp^*:H^2(G\times \SO(2)/T)\to H^2(G\times \SO(2)/F)\big)$ is generated by $e$. To find $\ker(\wp^*)$ consider the spectral sequences $E_r^{a,b}$ from $\wp$ and $\bE_r^{a,b}$ from $\bwp$. Choose the standard product generators $u_i\in H^1(T\simeq S^1\times S^1\times S^1)$ for $i=1,2,3$ and $v_i\in H^2(G\times \SO(2)/T\simeq S^2\times S^2)$ for $i=1,2$. Then it is clear the map $\bd_2^{0,1}: \bE_2^{0,1} \to \bE_2^{2,0}$ is given by $\bd_2^{0,1}(u_1)=v_1$, $\bd_2^{0,1}(u_2)=v_2$ and $\bd_2^{0,1}(u_3)=0$, for the right choice of sign on $v_i$. Further, since $F=\set{(e^{-i\alpha},e^{-i\beta}, R\big(n(-q\alpha+p\beta)\big))}$ we see $\pi^*:H^1(T/F)\to H^1(T)$ takes a certain generator $u\in H^1(T/F)$ to $nqu_1-npu_2-u_3$. Then by the commutativity of the diagram of spectral sequences induced by \eqref{eulerdiagram}, we see $d_2^{0,1}:E_2^{0,1}\to E_2^{2,0}$ satisfies $d_2^{0,1}(u)=\bd_2^{0,1}(\pi^*(u))=nqv_1-npv_2$. Therefore our induced map is given by $\wp^*:H^2(G\times \SO(2)/T)\simeq E_2^{2,0} =\grp{v_1,v_2}_\ab\twoheadrightarrow \grp{v_1,v_2}_\ab/\grp{nqv_1-npv_2}\simeq E_3^{2,0}\simeq H^2(G\times \SO(2)/F)$. In particular $\ker(\wp^*)=\grp{nqv_1-npv_2}=\grp{e}$ so $e=\pm n(qv_1-pv_2)\in \grp{v_1,v_2}_\ab\simeq H^2(G/L)$.

In particular this Euler class can take any possible nontrivial value, up to sign, for the right choices of $p,q,n$, where $\gcd(p,q)=1$ and $n>0$. (The case $e=0$ obviously corresponds to the trivial bundle $S^2\times S^2\times S^2\to S^2\times S^2$.) A change in the sign of $e$ changes the orientation of $P$ as an $\SO(2)$ bundle but does not effect $M\simeq P\times_{\SO(2)} S^2$. Hence any nontrivial $S^2$ sphere bundle $M$ over $S^2\times S^2$ with structure group $\SO(2)$ can be obtained in this manner for the right choice of $p,q,n$. Note, of course, that $M\to S^2\times S^2$ still may incidentally be trivial for certain values of $p,q,n$.


\subsection{Actions of type $N^6_C$:}

$$S^3\times S^3 \,\, \supset \,\, T^2, \,\, S^3\times \Z_n \,\, \supset \,\, S^1 \times \Z_n.$$

\vspace{.5em}


Letting $L=S^3\times S^1$ we see $G/L\simeq S^2$. We also find that $M_L$ is $S^4$, where $S^3\times S^1$ acts on $S^4\subset\R^3\times \R^2$ via $(g,e^{i\theta})\star (x,y)=(\psi(g)x,R(n\theta)y)$. Here $\psi:S^3\to \SO(3)$ is the usual double cover and $R(\theta)$ is rotation by $\theta$. Therefore, the non-primitivity bundle is $S^4\to M\to S^2$. 

From Section \ref{sec:nonprim} we know the structure group for this non-primitivity bundle is linear and hence $M$ is an $S^4$ sphere bundle over $S^2$. It is well know that there are only two equivalence classes of such bundles \cite{Steenrod}. To see which type $M$ is, consider the principal $\SO(5)$ bundle $P$ associated to the sphere bundle $M$. With the notation of Section \ref{sec:nonprim}, we can take $J=\SO(5)$ and $\phi:L\to \SO(5):(g,e^{i\theta})\mapsto \diag\big(\psi(g),R(n\theta)\big)$, the map induced from the action of $L$ on $M_L$ described above. Then $P\simeq G\times_L \SO(5)$.

We now claim that $P$ is simply connected if and only if $n$ is odd. From the homotopy sequence of the bundle $L\to G\times \SO(5) \to G\times_L \SO(5)$, it is clear that $P\simeq G\times_L \SO(5)$ is simply connected if and only if $i_*:\pi_1(L)\to \pi_1(G\times \SO(5))$ is onto. The claim then follows because $\pi_1(L)$ is generated by the curve $\alpha(\theta)=(1,e^{i\theta})$, and $i(\alpha(\theta))= \big(\alpha(\theta)^{-1},\diag(\Id_3,R(n\theta))\big)$ which generates $\pi_1(G\times \SO(5))\simeq \pi_1(\SO(5))$ if and only if $n$ is odd. Since $\SO(5)\times S^2$ is the principal $\SO(5)$ bundle associated to the trivial bundle $\pr_2:S^4\times S^2 \to S^2$ it follows that $M$ is diffeomorphic to $S^4\times S^2$ if $n$ is even and $M$ is diffeomorphic to the nontrivial $S^4$ sphere bundle over $S^2$ if $n$ is odd.

\subsection{Actions of type $N^6_D$:}

$$S^3\times S^3 \,\, \supset \,\, T^2, \,\, S^3\times S^1 \,\, \supset \,\, \set{(e^{ip\theta},e^{i\theta})}.$$

\vspace{.5em}


Again consider the non-primitivity bundle with $L=S^3\times S^1$. We see that $G/L\simeq S^2$ and that $M_L\simeq \CP^2$ where $L\simeq \SU(2)\times S^1$ acts on $\CP^2$ via $\phi:\SU(2)\times S^1\to \U(3): (A,z)\mapsto \diag(1, z^pA)$. So the non-primitivity bundle takes the form $\CP^2\to M\to S^2$ and the effective structure group can be taken to be $\U(3)/\set{zI|z\in S^1}\simeq \SU(3)/\Z_3$. Since $\pi_1(\SU(3)/\Z_3)\simeq \Z_3$ there are precisely three equivalence class of such bundles \cite{Steenrod}.

To distinguish these bundles consider their corresponding principal $\SU(3)/\Z_3$ bundles $P$. From Section \ref{sec:nonprim} we know that $P\simeq G\times_L (\SU(3)/\Z_3)$. Considering the bundle $L\to G\times (\SU(3)/\Z_3) \to G\times_L (\SU(3)/\Z_3)$, we see that $\pi_1(L)$ is generated by the curve $\alpha(\theta)=(I,e^{i\theta})$ whose image in $G\times (\SU(3)/\Z_3)$ is $(\alpha(\theta)^{-1}, [\diag(1, e^{ip\theta}I)])$. So $P$ is simply connected if and only if this curve generates $\pi_1(G\times (\SU(3)/\Z_3))$ which happens if and only if $p\equiv 1$ or $2 \mod 3$. It is clear that the principal bundle associated to the trivial bundle $\pr_2:\CP^2\times S^2\to S^2$ is $\SU(3)/\Z_3\times S^2$ which is not simply connected. We will now show that the principal bundles associated to the two nontrivial bundles are simply connected. In fact we claim that these two nontrivial bundles are diffeomorphic to each other and so are their associated principal bundles.

We can build the two nontrivial $\CP^2$ bundles over $S^2$ as follows \cite{Steenrod}. Choose a non-contractible curve $\gamma:S^1\to \SU(3)/\Z_3$, considering $S^1\subset \C$ as the standard unit circle. Then the curves $\gamma(z)$ and $\bar \gamma(z):=\gamma(\bar z)$ represent the two nontrivial elements of $\pi_1(\SU(3)/\Z_3)$. Note $\bar z$ denotes the complex conjugate of $z$ so $\bar\gamma$ is $\gamma$ traced backwards. Take two copies of the standard unit disk $D^2_\pm\subset \C$. Let $M_1=\CP^2\times D^2_-\cup \CP^2\times D^2_+/(x,z)\sim(\gamma(z)x,z)$ for $z\in S^1$ and $M_2=\CP^2\times D^2_-\cup \CP^2\times D^2_+/(x,z)\sim(\bar\gamma(z)x,z)$ for $z\in S^1$. It is clear how to extend this to make $M_1$ and $M_2$ smooth. We see $M_1$ and $M_2$ represent the two equivalence classes of nontrivial $\CP^2$ bundles over $S^2$ with group $\SU(3)/\Z_3$. Further, the function $f:M_1\to M_2$ given by $f(x,w)=(x,\bar w)$ is an explicit diffeomorphism between these spaces. Similarly the corresponding principal bundles are diffeomorphic, and hence they must both be simply connected since at least one of them is from above.

In conclusion, we have shown that $M$ must be diffeomorphic to $\CP^2\times S^2$ if $p\equiv 0\mod 3$ and diffeomorphic to either nontrivial $\CP^2$ bundle over $S^2$ with group $\SU(3)/\Z_3$ if $p\equiv 1$ or $2\mod 3$.

\subsection{Actions of type $N^6_E$:}

$$S^3\times S^3 \,\, \supset \,\, S^3\times S^1, \,\, S^3\times S^1 \,\, \supset \,\, \set{(e^{ip\theta},e^{i\theta})}.$$

\vspace{.5em}

Let $L=S^3\times S^1$ and note $G/L\simeq S^2$. Further, the primitive subdiagram is $S^3\times S^1 \supset S^3\times S^1,  S^3\times S^1 \supset  \set{(e^{ip\theta},e^{i\theta})}$ which corresponds to the action on $S^4$ given by $\phi:S^3\times S^1\to \SO(5):(g,e^{i\theta})\mapsto \diag(1,\psi(g,e^{ip\theta}))$. Here $\psi:S^3\times S^3\to \SO(4)$ is the usual double cover. So the non-primitivity bundle is $S^4\to M\to S^2$.

As in the case of $N^6_C$, we can find the diffeomorphism type of $M$ by considering the principal $\SO(5)$ bundle $P$ associated to the sphere bundle $M$. Again $P\simeq G\times_L \SO(5)$ where the map $i:L\to G\times \SO(5)$ is given by $i(x)=(x^{-1},\phi(x))$. The group $\pi_1(L)$ is again generated by the curve $\alpha(\theta)=(1,e^{i\theta})$ but this time $i(\alpha(\theta))= (\alpha(\theta)^{-1},\diag(1,R(-p\theta),R(p\theta))$ is contractible so $i_*:\pi_1(L)\to \pi_1(G\times \SO(5))$ is not onto. As with the case $N^6_C$, it follows that $M$ is diffeomorphic to $S^4\times S^2$.


\subsection{Actions of type $N^6_F$:}\label{N^6_F}

$$\SU(3) \,\, \supset \,\, \S(\U(2)\U(1)), \,\, \S(\U(2)\U(1)) \,\, \supset \,\, \SU(2)\SU(1)\cdot \Z_n$$
where $\Z_n\subset\S(\SU(1)\U(1)\U(1))$.

\vspace{.5em}

In this case let $L=\S(\U(2)\U(1))$. Clearly $G/L\simeq \CP^2$ and $M_L\simeq S^2$ where $L\simeq \U(2)$ acts on $S^2$ via $\phi:\U(2)\to \SO(3):A\mapsto R(\det(A)^n)$. Here $R(z)\in \SO(3)$ denotes rotation by $\arg(z)$ about some fixed axis, for $z\in S^1$. So the non-primitivity bundle is $S^2\to M\to \CP^2$ with structure group $\SO(2)\subset \SO(3)$. 

Consider the principal $\SO(2)$ bundle $P=G\times_L \SO(2)\to G/L$ associated to the sphere bundle $M$ and the reduced structure group $\SO(2)$. As before $P$ is determined by its Euler class $e\in H^2(G/L)\simeq \Z$ which we will compute. This process will be very similar to the case of $N^6_B$ so we will skip many details and simply indicate the important changes.

We see $P=G\times \SO(2)/F$ were $F=\set{(\ell^{-1},\phi(\ell))|\ell\in L}$, as before. In this case again, we have a bundle diagram similar to \eqref{eulerdiagram} with $T=\U(2)\times \SO(2)$, but this time with the top row of \eqref{eulerdiagram} replaced by the bundle $T/J\to G\times \SO(2)/J \to G\times SO(2)/T$ where $J=\SU(2)\times 1$. The key map $\pi^*:H^1(T/F\simeq S^1)\to H^1(T/J\simeq S^1\times S^1)\simeq \Z^2$ takes a generator to $\pm(n,1)$. Taking the spectral sequence $E_r^{a,b}$ for our bundle $\wp:P\to G/L$, we find $d_2^{0,1}(u)=nv$ for some generators $u\in E_2^{0,1}$ and $v\in E_2^{2,0}\simeq H^2(G/L)$, using the same method as for $N^6_B$. As in that case we see that $e=\pm nv\in H^2(G/L)\simeq \Z$, which can take on any nontrivial value, up to sign. The case $e=0$ corresponds to the trivial bundle. Again, the sign on $e$ does not effect $M\simeq P\times_{\SO(2)} S^2$. Therefore we can obtain any nontrivial $S^2$ sphere bundle $M$ over $\CP^2$ with reduced structure group $\SO(2)$.

{\setlength{\tabcolsep}{0.2cm}
\renewcommand{\arraystretch}{1.302}
\stepcounter{equation}
\begin{table}[!ht]
\begin{center}
\begin{tabular}{|c|c|}
\hline
$N^6_A$ & $S^3\times T^2 \,\, \supset \,\, \set{(e^{i a_-\theta},e^{ib_-\theta},e^{ic_-\theta})}\cdot H, \,\, \set{(e^{i a_+\theta},e^{ib_+\theta},e^{ic_+\theta})}\cdot H \,\, \supset \,\, H$\\
 & where $K^-\ne K^+$, $H=H_-\cdot H_+$, $\gcd(b_\pm,c_\pm)=1$,\\
 & $a_\pm=rb_\pm+sc_\pm$, and $K^-_0\cap K^+_0\subset H$\\
\cline{2-2}
 & $M \simeq S^3\times S^3$\\
\hline
\hline
$N^6_B$ & $S^3\times S^3 \,\, \supset \,\, \set{(e^{i\theta},e^{i\phi})}, \,\, \set{(e^{i\theta},e^{i\phi})} \,\, \supset \,\, \set{(e^{ip\theta},e^{iq\theta})} \cdot \Z_n$\\
 & where $\gcd(p,q)=1$ and $\Z_n\cap\set{(e^{ip\theta},e^{iq\theta})}=1$\\
\cline{2-2}
 & $S^2\to M\to S^2\times S^2$ with structure group $\SO(2)$\\
 & $e_P=\pm n(q,-p)\in H^2(S^2\times S^2)\simeq \Z\oplus \Z$\\
\hline
\hline
$N^6_C$ & $S^3\times S^3 \,\, \supset \,\, T^2, \,\, S^3\times \Z_n \,\, \supset \,\, S^1 \times \Z_n$\\
\cline{2-2}
 & $S^4\to M\to S^2$ with structure group $\SO(5)$\\
 & bundle trivial if and only if $n$ even\\
\hline
\hline
$N^6_D$ & $S^3\times S^3 \,\, \supset \,\, T^2, \,\, S^3\times S^1 \,\, \supset \,\, \set{(e^{ip\theta},e^{i\theta})}$\\
\cline{2-2}
 & $\CP^2\to M\to S^2$ with structure group $\SU(3)/\Z_3$\\
 & bundle trivial if and only if $p\equiv 0 \mod 3$\\
\hline
\hline
$N^6_E$ & $S^3\times S^3 \,\, \supset \,\, S^3\times S^1, \,\, S^3\times S^1 \,\, \supset \,\, \set{(e^{ip\theta},e^{i\theta})}$\\
\cline{2-2}
 & $M\simeq S^4\times S^2$\\
\hline
\hline
$N^6_F$ & $\SU(3) \,\, \supset \,\, \S(\U(2)\U(1)), \,\, \S(\U(2)\U(1)) \,\, \supset \,\, \SU(2)\SU(1)\cdot \Z_n$\\
 & where $\Z_n\subset\S(\SU(1)\U(1)\U(1))$\\
\cline{2-2}
 & $S^2\to M\to \CP^2$ with structure group $\SO(2)$\\
 & $e_P=\pm n\in H^2(\CP^2)\simeq \Z$\\
\hline
\end{tabular}
\end{center}
\vspace{0.1cm}
\caption{ {\bf Summary of results.} This table lists the group diagrams which describe a manifold $M$ of the given type, along with the information which determines the diffeomorphism type of $M$. For the $S^2$ bundles over $S^2\times S^2$ or $\CP^2$, $e_P$ denotes the Euler class of the associated principal $\SO(2)$ bundle. See Sections \ref{N^6_A} through \ref{N^6_F} for details.}\label{t:summary}
\end{table}}

%
%
%
%
%
%

\end{document}